\documentclass[journal]{IEEEtran}
\usepackage{graphicx}
\usepackage[normalem]{ulem}
\usepackage{amssymb}
\usepackage{amsmath}
\usepackage{amstext}
\usepackage{amsfonts}
\usepackage[latin1]{inputenc}
\usepackage[normalem]{ulem}
\usepackage[brazil]{babel}

\newtheorem{theorem}{Theorem}
\newtheorem{proposition}{Proposition}
\newtheorem{lemma}{Lemma}
\newtheorem{corollary}{Corollary}

\begin{document}

\title{Compactly Supported One-cyclic Wavelets Derived from Beta
Distributions}

\author{H. M. de Oliveira, G. A. A. de Ara\'{u}jo}

\maketitle

\begin{abstract}

\textbf{New continuous wavelets of compact support are introduced,
which are related to the beta distribution. They can be built from
probability distributions using "blur" derivatives. These new
wavelets have just one cycle, so they are termed unicycle wavelets.
They can be viewed as a soft variety of Haar wavelets whose shape is
fine-tuned by two parameters $a$  and $b$. Close expressions for
beta wavelets and scale functions as well as their spectra are
derived. Their importance is due to the Central Limit Theorem
applied for compactly supported signals.} \footnote{This work was
partially supported by the Brazilian National Council for Scientific
and Technological Development (CNPq) under research grant
$\#306180$. Authors are with Department of Electronics and Systems,
Federal University of Pernambuco, CTG-UFPE C.P 7800, 50.711-970,
Recife, Brazil, +55-81-2126-8210 fax: +55-81-2126-8215 (e-mail:
hmo@ufpe.br, giovanna.angelis@ibest.com.br).}
\end{abstract}


\IEEEpeerreviewmaketitle

\section{PRELIMINARIES AND BACKGROUND}
\PARstart{W}avelets are strongly connected with probability
distributions. Recently, a new insight into wavelets was presented,
which applies Max Born reading for the wave-function \cite{Beiser}
in such a way that an information theory focus has been achieved
\cite{de Oliveira2}. Many continuous wavelets are derived from a
probability density (e.g. Sombrero). This approach also sets up a
link among probability densities, wavelets and ``blur derivatives''
\cite{Kaiser}. To begin with, let $P(.)$ be a probability density,
$P\in \mathbb{C}^{\infty}$, the space of complex signals
$f:\mathbb{R}\rightarrow \mathbb{C}$ infinitely differentiable.

If \begin{equation}\stackrel \lim {_{t\rightarrow \infty}}
\frac{d^{n-1}P(t)}{dt^{n-1}}=0\end{equation} then

\begin{equation}\psi(t)= (-1)^n \frac{d^n P(t)}{dt^n}\end{equation}
is a wavelet engendered by $P(.)$. Given a mother wavelet $\psi$
that holds the admissibility condition [4,5] then the continuous
wavelet transform is defined by

\begin{equation} CWT(a,b) =  \int^{+\infty}_{-\infty} f(t) \cdot \frac{1}{\sqrt{|a|}}\psi(\frac{t-
b}{a})dt,\end{equation}

$\forall a \in \mathbb{R}-\{0\}$, $b\in \mathbb{R}$.

Continuous wavelets have often unbounded support, such as Morlet,
Meyer, Mathieu, de Oliveira wavelets [6-8]. In the case where the
wavelet was generated from a probability density, one has

\begin{equation}\frac{1}{\sqrt{|a|}}\psi_n(\frac{t-b}{a})=(-1)^n
\frac{1}{\sqrt{|a|}}\frac{\partial^n P(\frac{t-b}{a})}{\partial
t^n}.\end{equation}
 Now

\begin{equation}\frac{\partial^n P(\frac{t-b}{a})}{\partial b^n}= (-1)^n
\frac{1}{a^n} P^{(n)}(\frac{t-b}{a}),\end{equation}
  so that
\begin{equation}CWT(a,b)=\frac{1}{\sqrt{|a|}}
\int^{+\infty}_{-\infty} f(t) \cdot \frac{\partial^n P(\frac{t-
b}{a})}{\partial b^n}dt.\end{equation}

If the order of the integral and derivative can be commuted, it
follows that

\begin{equation} CWT(a,b)=
\frac{1}{\sqrt{|a|}} \frac{\partial^n }{\partial
b^n}\int^{+\infty}_{-\infty} f(t) \cdot P(\frac{t-
b}{a})dt.\end{equation}

Defining the LPFed signal as the "blur" signal

\begin{equation}\widetilde{f}(a,b)=\int^{+\infty}_{-\infty} f(t) \cdot \frac{1}{\sqrt{|a|}}P(\frac{t-
b}{a})dt=\int^{+\infty}_{-\infty} f(t) \cdot
P_{a,b}(t)dt,\end{equation}

an interesting interpretation can be made: set a scale $a$ and take
the average (smoothed) version of the original signal - the blur
version $\widetilde{f}(a,b)$. The ``blur derivative''

\begin{equation}\frac{\partial^n}{\partial
b^n}\widetilde{f}(a,b)\end{equation}

is the $n^{th}$ derivative regarding the shift $b$ of the blur
signal at the scale $a$. The blur derivative coincide with the
wavelet transform $CWT(a,b)$ at the corresponding scale. Details
(high-frequency) are provided by the derivative of the low-pass
(blur) version of the original signal.

\subsection{Revisiting Central Limit Theorems}

There are essentially three kinds of central limit theorems: for
unbounded distributions, for causal distributions and for compactly
supported distributions \cite{Gnedenko}. The random variable
corresponding to the sum of $N$ independent and identically
distributed (i.i.d.) variables converges to: a Gaussian
distribution, a Chi-square distribution or a Beta distribution (see
Table I). The Gaussian pulse always has been playing a very central
role in Engineering and it is associated with Morlet's wavelet,
which is known to be of unbounded support. This is the only wavelet
that meets the lower bound of Gabor's uncertainty inequality
\cite{Gabor}. The Gabor concept of logon naturally leads to the
Gaussian waveforms as an efficient signalling in the time-frequency
plan. Nevertheless, in the cases where a constraint in signal
duration is imposed, it can be expected that beta waveforms will be
the most efficient signalling in the time-frequency plan. The
concept of wavelet entropy was recently introduced and Morlet
wavelet also revealed to be a special wavelet [2,11]. Among all
wavelets of compact support, it can be expected that the one linked
to the beta distribution could also play a valuable practical and
theoretical role.

\begin{table}[!h]
\caption{DIFFERENT VERSIONS OF THE CENTRAL LIMIT THEOREM: UNBOUNDED
DISTRIBUTIONS, CAUSAL DISTRIBUTIONS AND COMPACTLY SUPPORTED
PROBABILITY DISTRIBUTIONS.} \label{tab1}
\begin{tabular}{c c}
\hline Marginal  Distribution & Central Limit Distribution as
$N\rightarrow\infty$\\ \hline
Unbounded Support&$G(t|m,\sigma^2)=\frac{1}{sqrt{2\pi\sigma^2}}\cdot e^{\frac{-(t-m)^2}{2\sigma^2}}$\\
Causal Distribution&$\chi^2(t|m,\sigma^2)=\frac{t^{\alpha}\cdot\beta^{-t\diagup\beta}}{\beta^{\alpha+1}\Gamma(\alpha+1)}$\\
Compact Support&$beta(t|\alpha,\beta)=K\cdot t^{\alpha} \cdot
(1-t)^{\beta},\quad0<t<1$\\ \hline
\end{tabular}
\end{table}

Let $p_i(t)$ be a probability density of the random variable
$t_i$,where $i=1,2,3..N$ i.e. $p_i(t)\geqslant0$  , $(\forall t)$
and

\begin{equation} \int^{+\infty}_{-\infty} p_i(t) dt=1.\end{equation}

If $p_i(t)\leftrightarrow P_i(\omega)$ , then $P_i(0)=1$  and
$(\forall \omega) |P_i(\omega)| \leq 1$ . Suppose that all variables
are independent. The density $p(t)$ of the random variable
corresponding
to the sum \begin{equation}t= \sum^N_{i=1}t_i \label{eleven}\end{equation} is given by the iterate convolution \cite{Davenport}\\
\begin{equation}p(t)=p_1(t)*p_2(t)*p_3(t)*...p_N(t). \end{equation}

If $p_i(t)\leftrightarrow P_i(\omega)=
|P_i(\omega)|e^{j\Theta_i(\omega)}$,  $i=1,2,3..N$ and
$p(t)\leftrightarrow P(\omega)= |P(\omega)|e^{j\Theta(\omega)}$ ,
then

\begin{equation}|P(\omega)|= \prod^N_{i=1}| P_i(\omega)|,\qquad\Theta(\omega)= \sum^N_{i=1}
\Theta_i(\omega).\end{equation}

The mean and the variance of a given random variable $t_i$ are,
respectively

\begin{equation} m_i =
 \int^{+\infty}_{-\infty} \tau \cdot p_i(\tau)d\tau,\end{equation}

\begin{equation} \sigma^2_i =
 \int^{+\infty}_{-\infty} (\tau-m_i)^2 \cdot p_i(\tau)d\tau.\end{equation}

The following theorems can be proved \cite{Gnedenko}.

\begin{theorem} \textit{Central Limit Theorem for distributions of unbounded
support}\\If the distributions $\{p_i(t)\}$  are not a lattice (a
Dirac comb) and $E(t^3_i)< \infty$, and
$$\lim_{N\rightarrow\infty}\sigma^2=+\infty$$ then $ t=\sum^N_{i=1}t_i$ holds, as $N \rightarrow \infty$,

\begin{equation}
P(\omega)\sim
e^{-{\frac{\sigma^2\omega^2}{2}}-jm\omega},\end{equation}

\begin{equation}p(t)\sim \frac {1}{\sqrt{2\pi\sigma^2}}
 e^{-{\frac{(t-m)^2}{2\sigma^2}}} \end{equation} $\Box$

\end{theorem}

According to Gnedenko and Kolmogorov, if all marginal probability
densities have bounded support, then the corresponding theorem is
\cite{Gnedenko}:

\begin{theorem} \textit{Central Limit Theorem for distributions of compact
support}. Let $\{p_i(t)\}$  be distributions such that $Supp\{(
p_i(t))\}=(a_i,b_i)(\forall i)$. Let

\begin{eqnarray}
a=\sum^N_{i=1}a_i<+\infty, \\ b=\sum^N_{i=1}b_i<+\infty.
\end{eqnarray}

It is assumed without loss of generality that $a=0$ and $b=1$. The
random variable defined by (\ref{eleven}) holds as
$N\rightarrow\infty$,

\begin{equation}
p(t)\sim \Biggl\lbrace \begin{array}{rcl}
k\cdot t^\alpha (1-t)^\beta ,&  0\leq t\leq 1\\
0, & \mbox{otherwise}\end{array}
\end{equation}

where

\begin{equation}
\alpha = \frac{m(m-m^2-\sigma^2)}{\sigma^2},
\end{equation}and
\begin{equation}
\beta = \frac{(1-m)(\alpha+1)}{m}.
\end{equation} $\Box$
\end{theorem}
In spite of the fact that a general theory of deriving wavelets from
probability distributions is well known, the particular application
discussed in this paper search for discovering a link between the
most noteworthy compact support distribution and wavelets.

\section {$\beta$-WAVELETS: NEW COMPACTLY SUPPORTED WAVELETS}

The beta distribution is a continuous probability distribution
defined over the interval $0\leq t\leq1$ \cite{Davies}. It is
characterized by a couple of parameters, namely  $\alpha$ and
$\beta$, according to:

\begin{equation}
P(t)=
\frac{1}{B(\alpha,\beta)}t^{\alpha-1}\cdot(1-t)^{\beta-1},\quad
1\leq\alpha,\beta\leq+\infty.
\end{equation}

The normalizing factor is $B(\alpha,\beta)=\frac{\Gamma(\alpha)\cdot
\Gamma(\beta)}{\Gamma(\alpha+\beta)}$, where $\Gamma(\cdot)$ is the
generalized factorial function of Euler and $B(\cdot,\cdot)$  is the
Beta function \cite{Davies}.

The following parameters can be computed:
\begin{eqnarray}
&Supp{(P)}=[0,1]\\
&mean=\frac{\alpha}{\alpha+\beta}\\
&mode=\frac{\alpha-1}{\alpha+\beta-2}\\
&variance=
\sigma^2=\frac{\alpha\beta}{(\alpha+\beta)^2(\alpha+\beta+1)}\\
&characteristic function=M(\alpha,\alpha+\beta,j\nu),
\end{eqnarray}

where $M(\cdot,\cdot,\cdot)$ is the Kummer confluent hypergeometric
function \cite{Slater,Gradshteyn}. The $N^{th}$ moment of $P(\cdot)$
can be found using

\begin{eqnarray} moment(N) =
\int^{1}_{0} t^N \cdot p_i(t)dt=
\frac{B(\alpha+N,\beta)}{B(\alpha,\beta)}\nonumber\\=
\frac{B(\alpha+\beta,N)}{B(\alpha,N)}.
\end{eqnarray}

The derivative of the beta distribution can easily be found.

\begin{equation}
P'(t)= \Bigl(\frac{\alpha-1}{t}-\frac{\beta-1}{1-t}\Bigl)P(t).
\end{equation}

A random variable transform can be made by an affine transform in
order to generate a new distribution with zero-mean and unity
variance \cite{Davenport}, which implies a non-normalized support
$T=\frac{1}{\sigma}=T(\alpha,\beta)$.\\Let a new random variable be
defined by $T\cdot(t-m)$. This variable has zero-mean and unity
variance. Its corresponding probability density is given by

\begin{eqnarray}
P(t|\alpha,\beta)=\frac{1}{B(\alpha,\beta)T(\alpha,\beta)}\Bigl
(\frac{t+m(\alpha,\beta)T(\alpha,\beta)}{T(\alpha,\beta)}\Bigl)^{\alpha-1}\cdot\nonumber\\
\Bigl(1-\frac{t+m(\alpha,\beta)T(\alpha,\beta)}{T(\alpha,\beta)}\Bigl)^{\beta-1}.
\end{eqnarray}

The $\beta$-wavelets can now be derived from these modified
distributions by using the concept of "blur" derivative. The
(unimodal) scale function associated with the wavelets is given by

\begin{equation}
\phi_{beta}(t|\alpha,\beta)=
\frac{1}{B(\alpha,\beta)T^{\alpha+\beta-1}}\cdot
(t-a)^{\alpha-1}\cdot (b-t)^{\beta-1},
\end{equation}
$a\leq t\leq b$. Since $P(\cdot|\alpha,\beta)$ is unimodal, the
wavelet generated by
\begin{equation}
\psi_{beta}(t|\alpha,\beta)=(-1) \frac{dP(t|\alpha,\beta)}{dt}
\end{equation}

has only one-cycle (a negative half-cycle and a positive
half-cycle).\\A close expression for first-order beta wavelets can
easily be derived. Within the support of $\psi(t|\alpha,\beta),\quad
a \leq t\leq b$,

\begin{eqnarray}
\psi_{beta}(t|\alpha,\beta)=\frac{-1}{B(\alpha,\beta)T^{\alpha+\beta-1}}\cdot
\Bigl[\frac{\alpha-1}{t-a}-\frac{\beta-1}{b-t}\Bigl]\cdot\nonumber\\(t-a)^{\alpha-1}\cdot
(b-t)^{\beta-1}
\end{eqnarray}
As a particular case, symmetric beta wavelets are given by
\begin{equation}
\psi_{beta}(t|\alpha,\alpha)= K(\alpha)\cdot t\cdot
[t^2-(2\alpha+1)]^{\alpha-2},
\end{equation}

where

\begin{equation}
K(\alpha)= (-1)^\alpha \cdot \frac
{2(\alpha-1)}{(2\sqrt{2\alpha+1})^{2\alpha-1}}\cdot \frac
{\Gamma(2\alpha)}{[\Gamma(\alpha)]^2}.
\end{equation}

The main features of beta wavelets of parameters $\alpha$ and
$\beta$ are:

\begin{eqnarray}
Supp (\psi) &=&\Bigl [ \frac{-1}{\sqrt{{\beta}\diagup{\alpha}}}
\sqrt{\alpha+\beta+1},\sqrt{\frac{\beta}{\alpha}}
\sqrt{\alpha+\beta+1}\Bigl ]\nonumber \\  &=& [a,b]
\end{eqnarray}

\begin{equation}
length Supp (\psi) = T(\alpha,\beta)=(\alpha+\beta)
\sqrt{\frac{\alpha+\beta+1}{\alpha\beta}}.
\end{equation}

The parameter $R= b\diagup |a|= \beta \diagup \alpha$  is referred
to as "cyclic balance", and is defined as the ratio between the
lengths of the causal and non-causal piece of the wavelet. It can be
easily shown that the instant of transition $t_{zerocross}$ from the
first to the second half cycle is given by

\begin{equation}t_{zerocross}= \frac{(\alpha-\beta)}{(\alpha+\beta-2)}
\sqrt{\frac{\alpha+\beta+1}{\alpha\beta}}.
\end{equation}

Although scale and wavelets can be found for any $\alpha,\beta >1$,
the behavior of the wavelet at the extreme points of the support can
be discontinuous (e.g. see Table Ia). However, it is a simple matter
to guarantee the continuity of the wavelet according to:

\begin{proposition}
Beta one-cycle wavelets of parameters $\alpha>2$ and  $\beta>2$ are
smooth, continuous wavelets of compact support.\\
\begin{proof}Clearly, $\psi_{beta}(t|\alpha,\beta)=0\quad\forall
t<a$ and $\forall t>b$. The only concerns are therefore with the
extreme points of the support, but
$\psi_{beta}(a|\alpha,\beta)=\psi_{beta}(b|\alpha,\beta)=0$ provided
that $\alpha>2$ and $\beta>2.$
\end{proof}
\end{proposition}

Remember that $ 1 \leq \alpha,\beta \leq +\infty $ and $a<0<b$.\\The
beta wavelet spectrum can be derived in terms of the Kummer
hypergeometric function \cite{Slater}, which is solution of the
equation

\begin{equation}
z\frac{d^2\omega}{dz}+(\beta-z)\frac{d\omega}{dz}-\alpha\omega=0.
\end{equation}

Let $\psi_{beta}(t|\alpha,\beta)\leftrightarrow
\Psi_{BETA}(\omega|\alpha,\beta)$ denote the Fourier transform pair
associated with the wavelet. This spectrum is also denoted by
$\Psi_{BETA}(\omega)$ for short. It can be proved by applying
properties of the Fourier transform that

\begin{eqnarray}
\Psi_{BETA}(\omega)&=-j\omega \cdot M \Bigl(
\alpha,\alpha+\beta,-j\omega
(\alpha+\beta)\sqrt{\frac{\alpha+\beta+1}{\alpha\beta}}\Bigl)\cdot\nonumber\\
&exp{(j\omega \sqrt{\frac{\alpha(\alpha+\beta+1)}{\beta}})}.
\end{eqnarray}

The spectrum of a number of unicyclic beta wavelets is presented in
Figure 1. The spectrum evaluation was carried out using the
relationship:

\begin{equation}
M(\alpha,\alpha+\beta,j\nu)=
\frac{\Gamma(\alpha+\beta)}{\Gamma(\alpha) \cdot \Gamma(\beta)}\cdot
\int^1_0 e^{j\nu t}t^{\alpha-1}(1-t)^{\beta-1}dt
\end{equation}

Only symmetrical $(\alpha=\beta)$ cases have zeroes in the spectrum
(Fig. 1a). A few asymmetric $(\alpha \neq \beta)$ beta wavelets are
shown in Fig. 1b. Inquisitively, they are parameter-symmetrical in
the sense that they hold
\begin{equation}|\Psi_{BETA}(\omega|\alpha,\beta)|=|\Psi_{BETA}(\omega|\beta,\alpha)|.\end{equation}
The spectrum of symmetric beta wavelets has been compared to that of
Haar wavelets of same support, just to check on the reliability of
spectral computations. The first spectral null of a Haar wavelet of
support $T(\alpha,\beta)$ occurs at a frequency $\omega_0=
\frac{2\pi}{T(\alpha,\beta)\diagup2}$ or $\nu_0=4\pi$, then at $6\pi
, 8\pi$ etc. For $\alpha =3$, the first spectral null occurred at a
frequency $\nu=11.526918406...$, which is close to $\nu_0=4\pi$ as
expected. As $\alpha$ increases, the wavelet half-cycle tends to be
shrunk (e.g. fig.2 b and e), thereby increasing the frequency of the
first spectral notch (Fig.1a).

\begin{figure}[!h]
\centering
\includegraphics[width=0.5\textwidth]{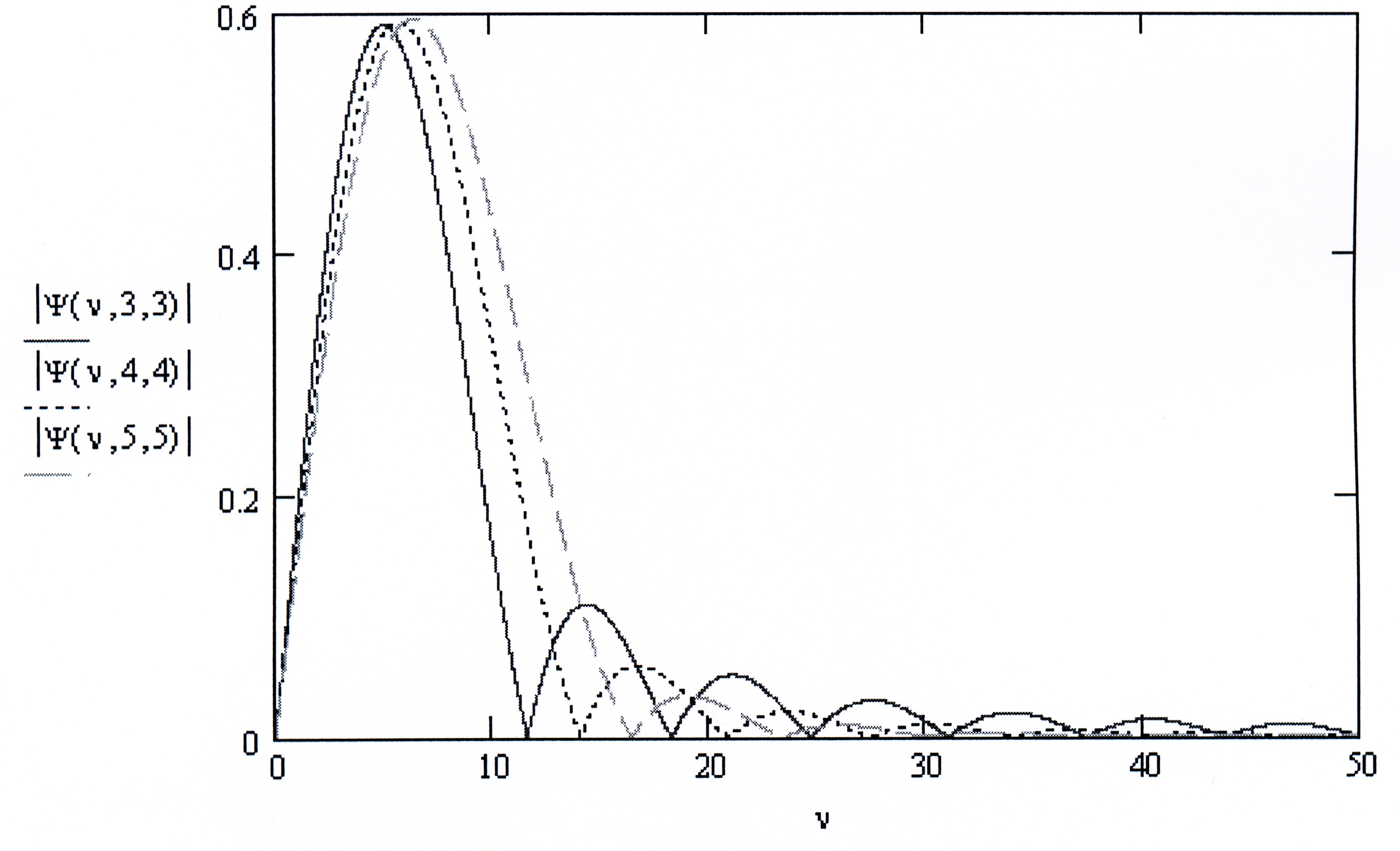}
\includegraphics[width=0.5\textwidth]{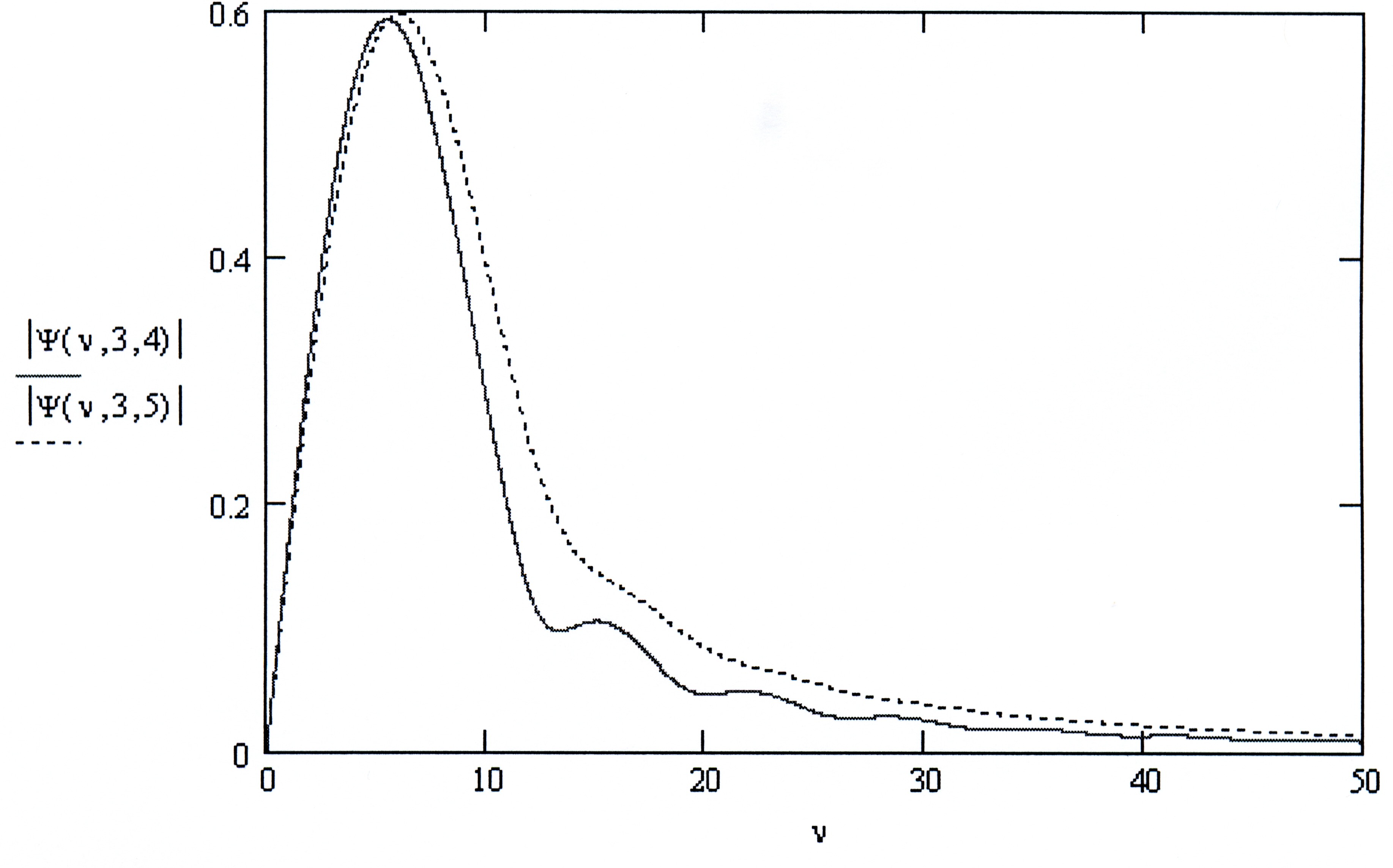}
\caption{ Magnitude of the spectrum  $\Psi_{BETA}(\omega)$ of a few
beta wavelets, $|\Psi_{BETA}(\nu|\alpha,\beta)|\times \nu$ for: a)
symmetric beta wavelets $\alpha = \beta=3$ (solid),  $\alpha =
\beta=4 $(dot) and $\alpha = \beta=5$ (dash); b) asymmetric beta
wavelets $\alpha=3$, $\beta=4$ (solid) and $\alpha=3$, $\beta=5$
(dot). In both plots the frequency axis is a normalized spectral
frequency $\nu=\omega T(\alpha,\beta)$.}
\end{figure}

\begin{figure}[!h]
\centering
\includegraphics[width=0.5\textwidth]{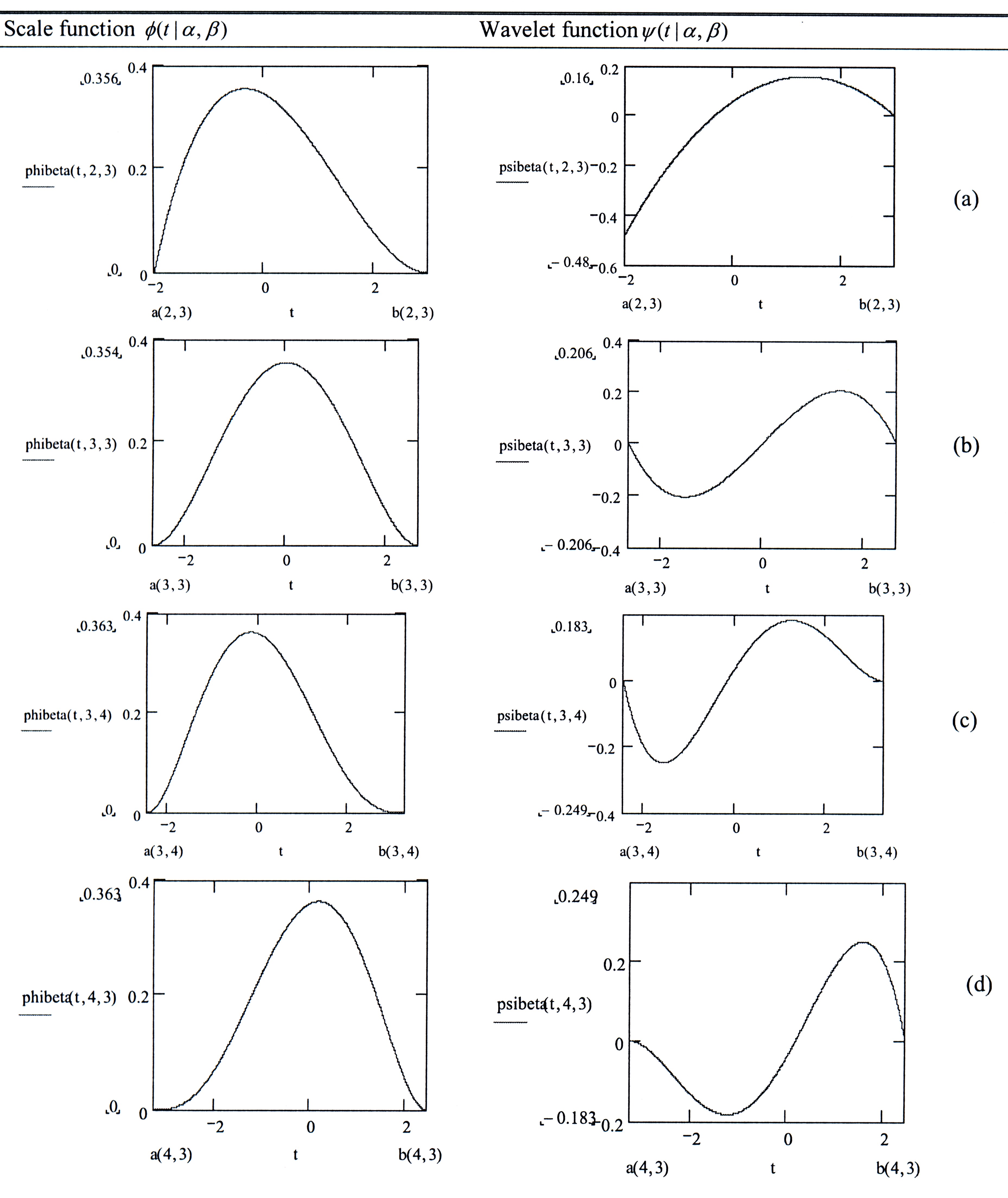}
\includegraphics[width=0.5\textwidth]{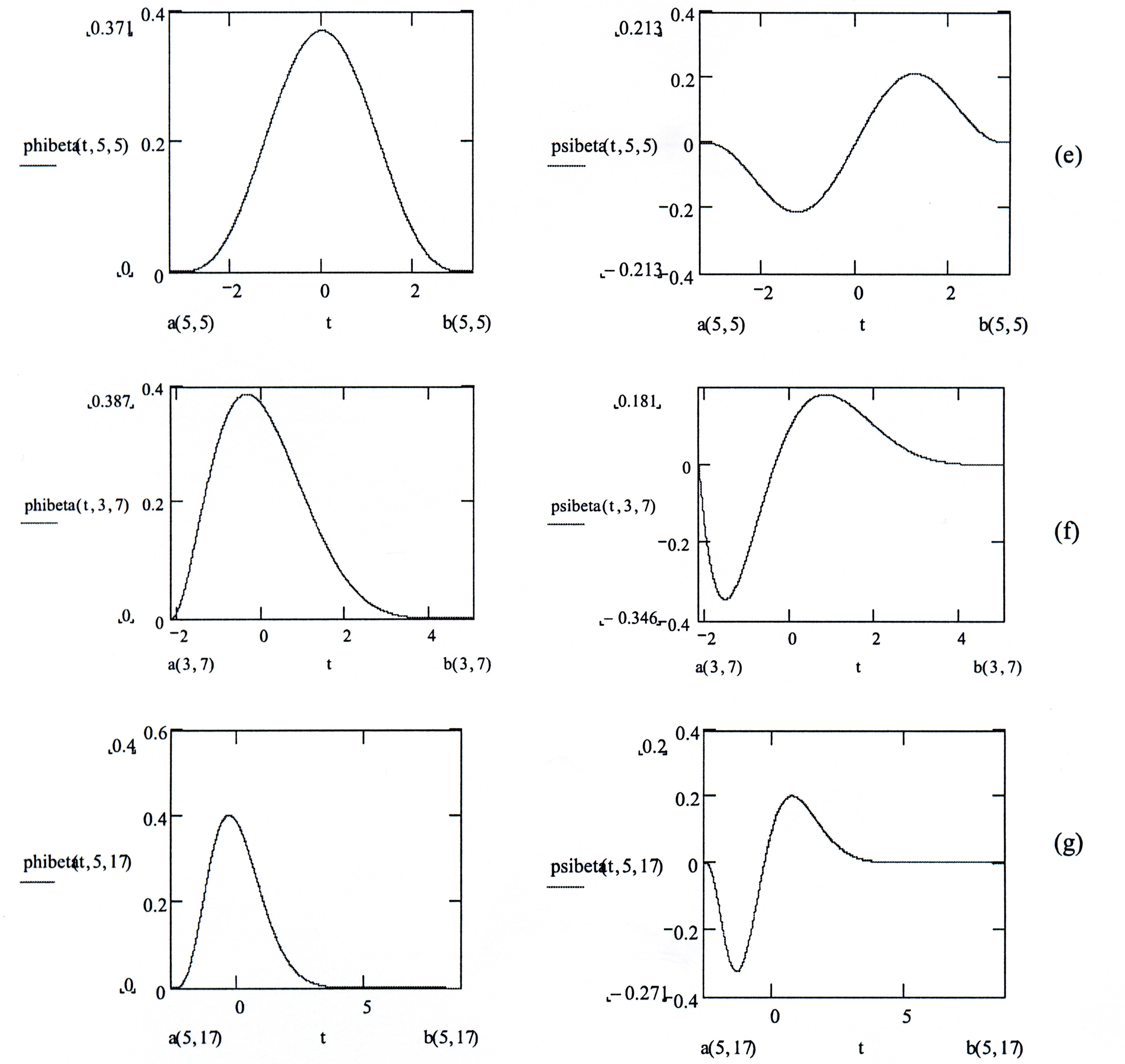}
\caption {Unicycle-beta scale function and wavelet for different
parameters:a) $\alpha=2,\beta=3$ b) $\alpha=\beta=3$ c) $\alpha=3,
\beta=4$ d) $\alpha=4,\beta=3$ e) $\alpha=\beta=5$ f)
$\alpha=3,\beta=7$ g) $\alpha=5,\beta=17$.} \label{fig2}
\end{figure}

\section{HIGH-ORDER BETA WAVELETS}

Due to the unimodal feature of the beta distribution, its first
derivative has just one cycle. Higher derivatives may also generate
further beta wavelets. Higher order beta wavelets are defined by

\begin{equation}
\psi_{beta}(t|\alpha,\beta)=(-1)^N
\frac{d^NP(t|\alpha,\beta)}{dt^N}.
\end{equation}

This is henceforth referred to as an $N$-order beta wavelet. They
exist for order $N\leq Min(\alpha,\beta)-1$. After some algebraic
handling, their close expression can be found:

\begin{eqnarray}
\Psi_{beta}(t|\alpha,\beta)&=&\frac{(-1)^N}{B(\alpha,\beta)\cdot
T^{\alpha+\beta-1}} \sum^N_{n=0} sgn(2n-N)\cdot\nonumber\\&
&\frac{\Gamma(\alpha)}{\Gamma(\alpha-(N-n))}
(t-a)^{\alpha-1-(N-n)}\cdot\nonumber\\& &
\frac{\Gamma(\beta)}{\Gamma(\beta-n)}(b-t)^{\beta-1-n}.
\end{eqnarray}
The choice of the order $N$ plays some role in the regularity of the
beta wavelets, and might be related with the Hölder and Sobolev
regularity. This topic, however, is not addressed in this paper. \\
A couple of high beta wavelets are shown in Fig.3. With the aim of
allowing the investigation of some potential applications of such
wavelets, routines to compute them should be written.

\begin{figure}[!htb]
\centering
\includegraphics[width=0.5\textwidth]{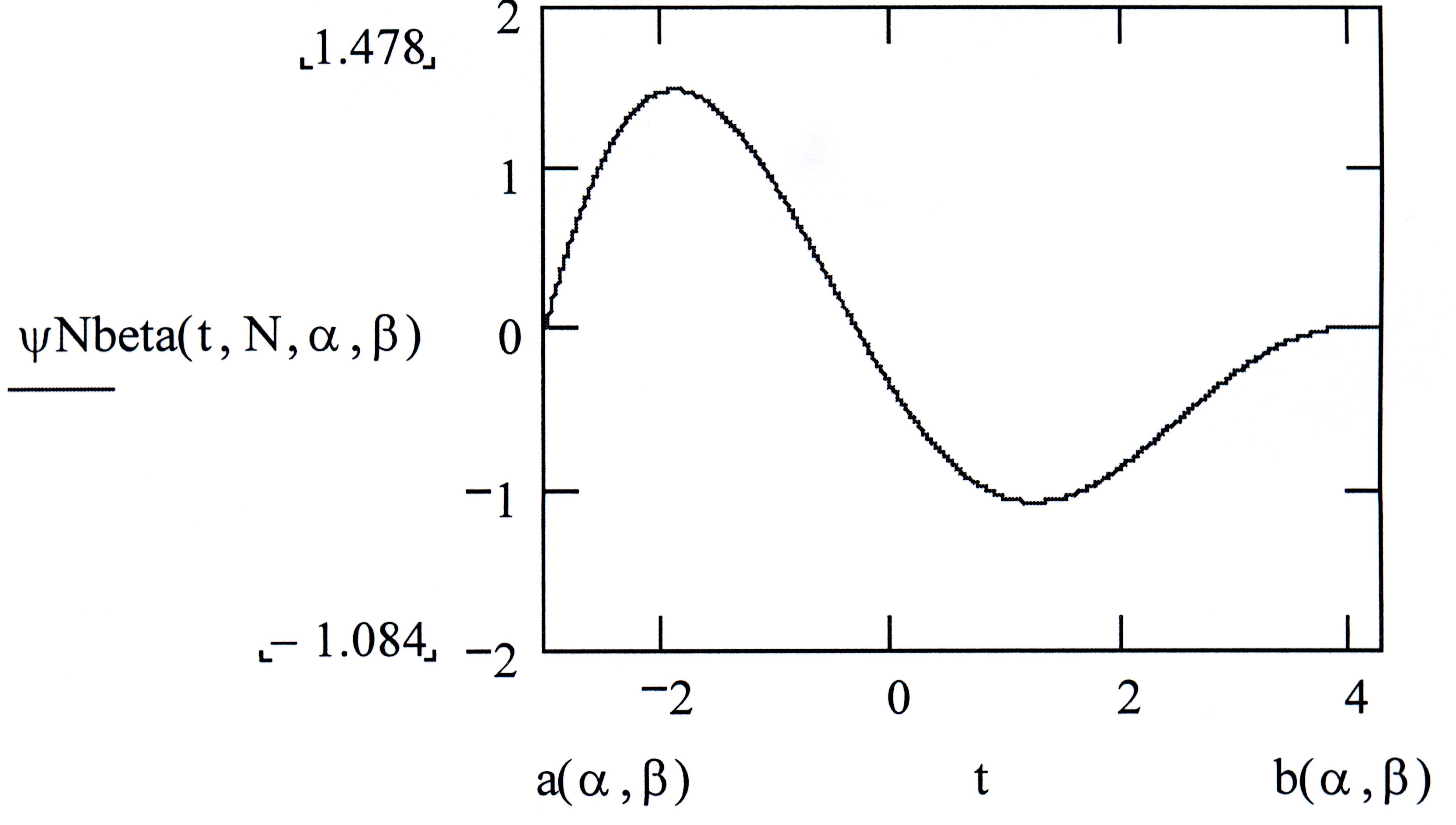}
\includegraphics[width=0.5\textwidth]{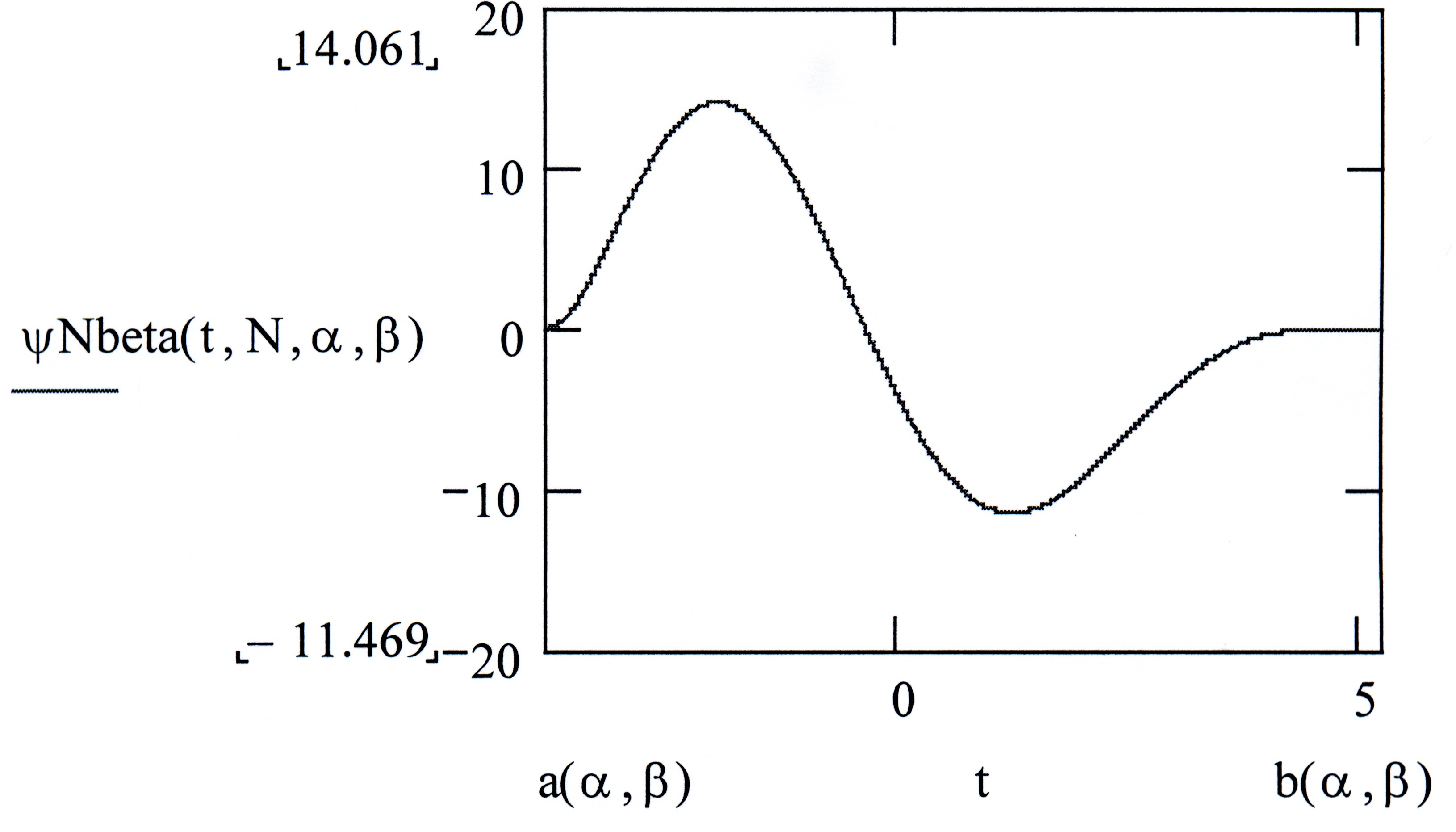}
\caption{High-order beta wavelets for different parameters: a)
$N=3$, $\alpha=5,\beta=7$; b) $N=5$, $\alpha=8,\beta=11$.}
\end{figure}
As it happens with any wavelet that has compact support, beta
wavelets share the ability of providing good estimate for short
transient, because no matter how short the interval is, there is a
scaled wavelet version whose support is limited within this time.
They can thus model local features efficiently as they are not
concerned with the data behavior far way from the focused location.
The main drawback of Haar wavelets is their discontinuities, which
engender a broadband spectrum. In contrast, beta wavelets can
provide a better balance between time and frequency resolution due
to their soft shape. For a given support length (time resolution), a
beta wavelet provides narrow spectrum (frequency resolution) than
the corresponding Haar wavelet of the same support. Nowadays one of
the most powerful software supporting wavelet analysis is the
$Matlab^{TM}$ \cite{Kamen}, especially when the wavelet graphic
interface is available. In the $Matlab^{TM}$ wavelet toolbox, there
exist five kinds of wavelets (type the command $waveinfo$ on the
prompt): (i) crude wavelets, (ii) infinitely regular wavelets, (iii)
orthogonal and compactly supported wavelets, (iv) biorthogonal and
compactly supported wavelet pairs, (v) complex wavelets. Figure 4
illustrates the beta wavelet implementation over
Matlab\texttrademark. The m-files that allow the computation of the
beta wavelet transform are currently $(freeware)$ available at the
URL: http://www2.ee.ufpe.br/codec/WEBLET.html (new wavelets).

\begin{figure}[!t]
\centering
\includegraphics[width=0.5\textwidth]{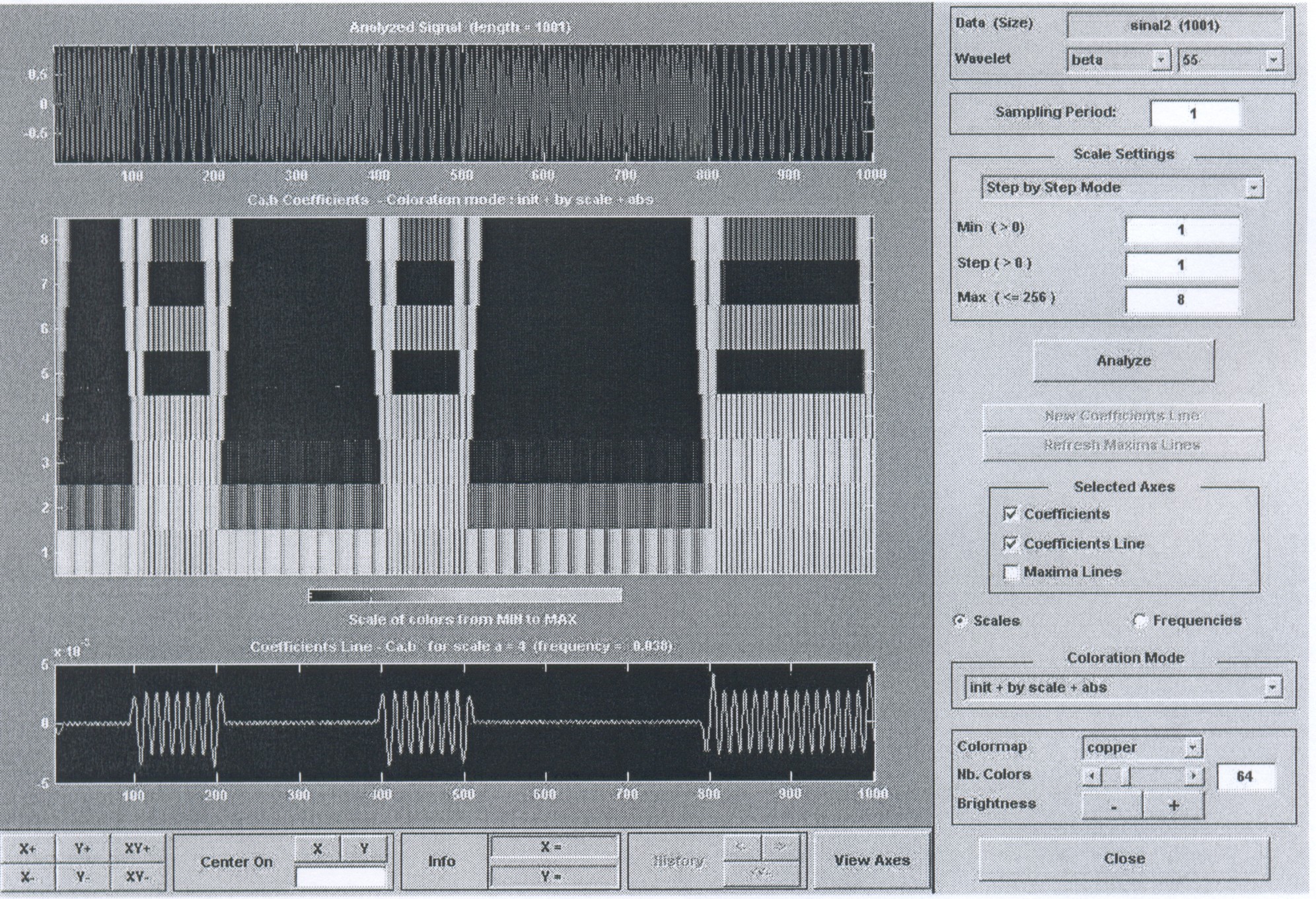}
\caption{Beta wavelets displayed over $Matlab^{TM}$ using the
\emph{wavemenu} command. The analyzed signal is a binary FSK signal.
This wavelet decomposition can be used to implement an efficient
frequency discriminator.}
\end{figure}

\section{CONCLUDING REMARKS}{Compactly supported wavelets are among
the most functional and useful wavelets. This correspondence
introduces a new family of wavelets of this class. These wavelets
can be viewed as some kind of soft-Haar wavelets. This new family of
wavelets looks like a sort of soft-Haar wavelets, since both are
unicycle wavelets. Besides the fact that beta wavelets are smoother,
they have extra flexibility since the balance between the two
half-cycles can be fine-tuned. However, it should be kept on mind
that such a comparison is rather loose, since completeness and
orthogonality properties that Haar wavelet holds were not addressed
for beta wavelets. It remains to be investigated how beta wavelets
can be approximated using FIR or IIR filters. In comparison with
other wavelets of compact support (e.g. dBN, coiflets etc.), the
beta wavelets derived in this work have
 idiosyncrasies and advantages: i) They are regular and smooth,
ii) have only one cycle, iii) have an analytical formulation (close
formulae), iv) Their importance rely on the Central Limit Theorem.
Many practical signals have a cyclic true nature, that is, they are
composed by uninterrupted cycles, as it occurs with power line
signals, some biomedical signals or FSK signals. Since many
waveforms are inborn-generated by successive cycles, their local
properties can probably better investigated via a wavelet able to
cope the changes from one cycle to another. It is often crucial to
detect a small discrepancy within a cycle. This behavior can be
useful for analyzing signals from certain modulation schemes or from
power systems disturbances. Particularly, Beta-wavelet-based FSK
modulation schemes are currently on investigation.}

\section{ACKNOWLEDGMENTS}
The authors thank Dr. Renato Jos\'e de Sobral Cintra for quite a lot
of comments, particularly regarding repeated convolution and Central
Limit Theorems.

\appendix
\begin{lemma} The square of a normalized beta density
\begin{equation}P(t)=\phi(t|\alpha,\beta)= \frac{1}{B(\alpha,\beta)}\cdot
t^{\alpha-1} (1-t)^{\beta-1} \end{equation}is proportional to
another beta density.
\end{lemma}

\begin{proof}
A straightforward algebraic handling yields
$\phi^2(t|\alpha,\beta)=\lambda_0 \cdot \phi(t|2\alpha-1,2\beta-1)$,
where
\begin{equation}
\lambda_0=\lambda_0(\alpha,\beta)=\frac{B(2\alpha-1,2\beta-1)}{B^2(\alpha,\beta)}
\end{equation}
\end{proof}

Let $D^{beta}=\{\phi_{beta}(t|\alpha,\beta)\}_{\alpha,\beta \in
\mathbb{R}}$ be the set of all possible signals of the kind beta
probability density (possibly non-normalized).
\begin{lemma}
The square of any beta density is proportional to another
beta-shaped density of same support.
\end{lemma}

\begin{proof}
The support of $\phi_{beta}(t|\alpha,\beta)$ remains unchanged and
furthermore,
$\phi^2_{beta}(t|\alpha,\beta)=\lambda\phi(t|2\alpha-1,2\beta-1),\quad
t \in [a,b]$, where
\begin{eqnarray}
\lambda_0=\lambda_0(\alpha,\beta) &=&
{\frac{T(2\alpha-1,2\beta-1)^{2(\alpha+\beta-1)}}{T^2(\alpha,\beta)\cdot
T(2\alpha-1,2\beta-1)}}.
\end{eqnarray}

\end{proof}

\begin{corollary}
$D^{beta}$ is a closed class of signals regarding the following
operations: rising to a power (pair exponent) and repeated
convolution (a number pair of times), i.e.
$\phi_{beta}(t|\alpha,\beta)\in D^{beta}\Rightarrow
\phi^2_{beta}(t|\alpha,\beta)\in D^{beta}$ and
$\phi_{beta}(t|\alpha,\beta)\ast\phi_{beta}(t|\alpha,\beta)\in
D^{beta}$.
\end{corollary}

A similar property is shared with the other densities concerned with
versions of the Central Limit Theorem.

\begin{lemma}
$\frac{1}{2\pi}\int^\infty_{-\infty}|M(\alpha,\alpha+\beta,j\nu)|^2d\nu=\lambda_0(\alpha,\beta).$
\end{lemma}

\begin{proof}
Parseval's identity furnishes
\begin{equation}\frac{1}{2\pi}\int^\infty_{-\infty}|M(\alpha,\alpha+\beta,j\nu)|^2d\nu=\int^1_0
\phi_{beta}^2(t|\alpha,\beta)dt\end{equation}and the proof follows
by applying lemma 1.
\end{proof}

\begin{lemma}
The second moment of the square of the Kummer hypergeometric
function $M(\alpha,\alpha+\beta,j\nu)$ is given by
\begin{equation}\frac{1}{2\pi}\int^\infty_{-\infty}\nu^2|M(\alpha,\alpha+\beta,j\nu)|^2d\nu=\chi(\alpha,\beta),
\end{equation}where
\begin{eqnarray}
\chi(\alpha,\beta)&=&\Bigl(\frac{1}{B^2(\alpha,\beta)}\Bigl)\Bigl[(\alpha-1)^2B(2\alpha-3,2\beta-3)\nonumber\\&
&-2(\alpha-1)(\alpha+\beta-2)B(2\alpha-3,2\beta-3)\nonumber\\&
&+(\alpha+\beta-2)^2B(2\alpha-3,2\beta-3)\Bigl].
\end{eqnarray}
\end{lemma}

\begin{proof}
It follows from Parseval's identity that
\begin{equation}\frac{1}{2\pi}\int^\infty_{-\infty}|\nu
M(\alpha,\alpha+\beta,j\nu)|^2d\nu=\int^1_0
\Bigl[\frac{d\phi(t|\alpha,\beta)}{dt}\Bigl]^2dt.\end{equation} Now
\begin{equation}
\frac{d\phi(t|\alpha,\beta)}{dt}= \frac{1}{B(\alpha,\beta)} \cdot
\Bigl[\frac{\alpha-1}{t}-\frac{\beta-1}{1-t}\Bigl]\cdot
t^{\alpha-1}\cdot (1-t)^{\beta-1}
\end{equation}and therefore the evaluation of the integral

\begin{equation}
\int^1_0\Bigl[\frac{(\alpha-1)-(\alpha+\beta-2)\cdot
t}{t(1-t)}\Bigl]^2\cdot t^{2\alpha-2}\cdot (1-t)^{2\beta-2}dt
\end{equation}completes the proof.
\end{proof}

The energy of the beta scale $\phi_{beta}(t|\alpha,\beta)$ and
wavelet function $\psi_{beta}(t|\alpha,\beta)$ can be computed
according to the following proposition.

\begin{proposition}The energies of the beta scale and wavelet function are, respectively,
\begin{equation}
E_{\phi_{beta}}=
\int^\infty_{-\infty}\phi^2_{beta}(t|\alpha,\beta)dt=\frac{\lambda_0(\alpha,\beta)}{T(\alpha,\beta)}
\end{equation} and \begin{equation}
E_{\psi_{beta}}=
\int^\infty_{-\infty}\psi^2_{beta}(t|\alpha,\beta)dt=\frac{\chi(\alpha,\beta)}{T^3(\alpha,\beta)}.
\end{equation}
\end{proposition}

\begin{proof}
A simple variable change gives \begin{equation} E_{\phi_{beta}}=
\int^{b}_{a}\phi^2_{beta}(t|\alpha,\beta)dt=\frac{1}{T(\alpha,\beta)}\int^{b}_{a}\phi^2(t|\alpha,\beta)dt
\end{equation} and the proof of the first part
follows by applying lemma 1. Let $\mathfrak{F}$ denote the Fourier
transform operator. Parseval's identity can be used in order to
evaluate $E_{\psi_{beta}}$:
\begin{equation}
\int^{b}_{a}\psi^2_{beta}(t|\alpha,\beta)dt= \frac{1}{2\pi}
\int^{\infty}_{-\infty}\Bigl|\mathfrak{F}\Bigl(\frac{dP(t|\alpha,\beta)}{dt}\Bigl)\Bigl|^2d\omega,
\end{equation}
so that
\begin{equation}
E_{\psi_{beta}}=\frac{1}{2\pi}\int^\infty_{-\infty}\Bigl|\omega
\cdot M(\alpha,\alpha+\beta,-j\omega \cdot
T(\alpha,\beta))\Bigl|^2d\omega.
\end{equation}
By a suitable variable change $\nu=\omega \cdot T(\alpha,\beta)$,
\begin{equation}E_{\psi_{beta}}=\frac{1}{2\pi T^3(\alpha,\beta)}\int^\infty_{-\infty}\nu^2\Bigl|M(\alpha,\alpha+\beta,j\nu)\Bigl|^2d\nu
\end{equation} and the proof follows from lemma 4.
\end{proof}

\begin{proposition}
Let $\alpha>1$ and $\beta>1$. The admissibility constant $c_{\Psi}$
of a unicyclic beta wavelet is
\begin{equation}
c_{\psi}(\alpha,\beta)=\frac{2\pi
\lambda_0(\alpha,\beta)}{T(\alpha,\beta)}<+\infty
\end{equation}
\end{proposition}
\begin{proof}
Since that
\begin{equation}c_{\psi}=\int^\infty_{-\infty}\frac{|\Psi_{BETA}(\omega)|^2}{|\omega|}d\omega,\end{equation}
then
\begin{equation}c_\Psi(\alpha,\beta)=\int^\infty_{-\infty}\Bigl|\omega
M(\alpha,\alpha+\beta,-j\omega \cdot T(\alpha,\beta))\Bigl|^2d\omega
\end{equation}the proof is completed using the lemma 3.
\end{proof}
Further interesting properties of beta distributions can be found in
\cite{Krysicki}.

\newpage
\begin{minipage}[t]{1\linewidth}
\begin{biography}[{\includegraphics[width=1in,height=1.25in,clip,keepaspectratio]{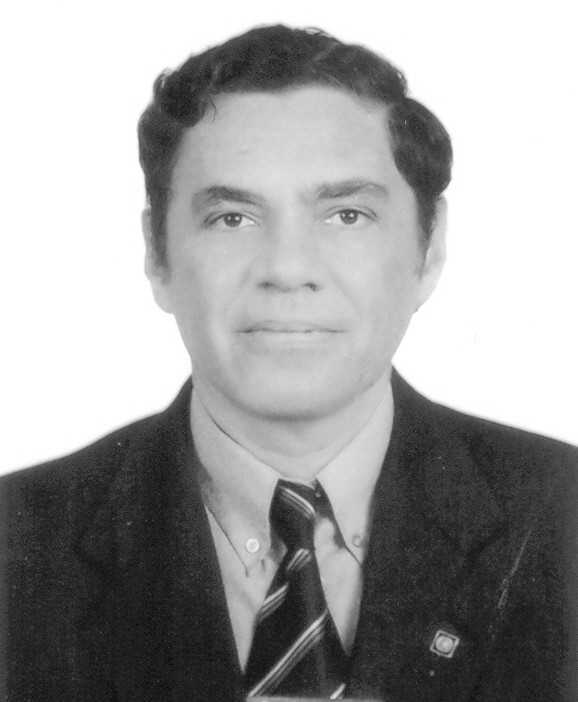}}]
{H\'elio Magalh\~{a}es de Oliveira} was born in Arcoverde,
Pernambuco, Brazil, in 1959. He received both the B.Sc. and M.S.E.E.
degrees in electrical engineering from Universidade Federal de
Pernambuco (UFPE), Recife, Pernambuco, in 1989 and 1983,
respectively. Then he joined the staff of the Department Electronics
And Systems (DES-UFPE) at the same university as a lecturer. In 1992
he earned the Docteur de l\emph{\'Ecole Nationale Sup\'erieure des
T\'el\'ecommunications} degree, in Paris. Dr. de Oliveira was
appointed as honored professor by twenty electrical engineering
undergratuate classes and chosen as the godfather of five
engineering graduation. His publications are available (in the .pdf
format) at http://www2.ee.ufpe.br/codec/publicacoes.html. He was the
head of the UFPE electrical engineering graduate program from 1992
to 1996. Research current interests include: communications theory,
applied information theory, genomic signal processing, signal
analysis, and wavelets. He was member of IEEE and is currently a
member SBrT (Brazilian Telecommunications Society). Dr. H\'elio de
Oliveira is currently an Division Editor of JBTS (J. of the
Brazilian Telecommunication Society).
\end{biography}
\begin{biography}[{\includegraphics [width=1in,height=1.25in,clip,keepaspectratio]{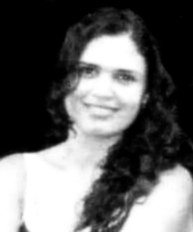}}]
{Giovanna Angelis Andrade de Ara\'ujo} Received the B.S.E.E. degree
from Universidade Federal da Para\'iba, Campina Grande,PB, Brazil,
in 1994. She is presently a M.S. student in electrical and
electronics engineering at the Universidade Federal de Pernambuco.
Her research interests include wavelets, fractals, and statistical
signal and image processing.
\end{biography}
\end{minipage}
\end{document}